# A forgotten little chapter on isoperimetric inequalities: On the fraction of a convex and closed plane area lying outside a circle with which it shares a diameter

José M. Pacheco
Departamento de Matemáticas[1]
Universidad de Las Palmas de Gran Canaria

**Abstract**

Often some interesting or simply curious points are left out when developing a theory. It seems that one of them is the existence of an upper bound for the fraction of area of a convex and closed plane area lying outside a circle with which it shares a diameter, a problem stemming from the theory of isoperimetric inequalities. In this paper such a bound is constructed and shown to be attained for a particular area. It is also shown that convexity is a necessary condition in order to avoid the whole area lying outside the circle.

**Key Words**: convex plane area, isoperimetric inequality.

**AMS Classification numbers**: 52A10, 52A38.

## Introduction

Possibly one of the oldest extremal problems is to find a set in Euclidean space with given surface area and enclosing maximum volume. There exist considerable differences in the mathematical treatment of the cases $n = 2$ and general $n$, as shown in the classical reference [2], where the general setting directly invites the reader to the realm of geometric measure theory. These problems pervade mathematical activity and many mathematicians have dealt with them to different depth degrees. As the motivating example for this paper, on reading the book *Littlewood's Miscellany* [1] one finds in page 32 the following observation:

"*An isoperimetrical problem*: an area of (greatest) diameter not greater than 1 is at most $\frac{1}{4}\pi$"

Littlewood's "greatest diameter" is now better known as *the* diameter $d$ of the plane area, defined as

$$d = \sup\{dist(X,Y) | X, Y \in boundary\},$$

and without loss of generality, we can suppose in Littlewood's remark the area to be convex and bounded by a continuous closed curve. Indeed non-convexity would only amount to reducing the enclosed area, thus enhancing the inequality (see figure 1).

---

[1] Postal address: Campus de Tafira Baja, 35017 LAS PALMAS, Spain. E-mail: pacheco@dma.ulpgc.es

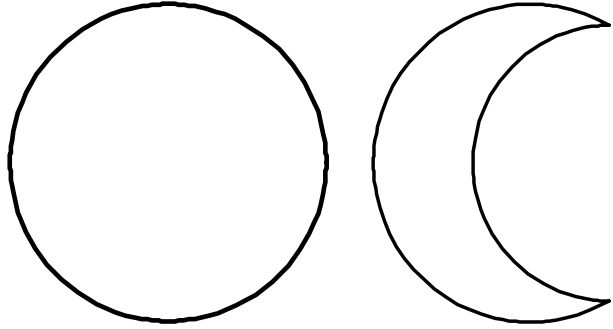

Figure 1: Nonconvexity implies area loss.

Littlewood's argument is the following (figure in p. 33 of [1]), see figure 2:

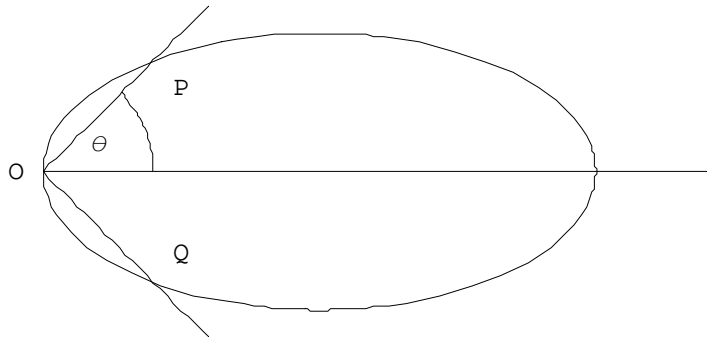

Figure 2: A recreation of the figure in [1], p. 33.

$$area = \frac{1}{2}\int_0^{\frac{\pi}{2}}(OP^2 + OQ^2)d\theta$$

where $OP = \rho(\theta)$, $OQ = \rho(\theta - \frac{\pi}{2})$, and $OP^2 + OQ^2 = PQ^2 \leq diam^2 \leq 1$. Therefore

$$area = \frac{1}{2}\int_0^{\frac{\pi}{2}} PQ^2 d\theta \leq \frac{1}{2}\int_0^{\frac{\pi}{2}} d\theta = \frac{\pi}{4}$$

The bound $\frac{1}{4}\pi$ is attained for the unit diameter circle, and the classical isoperimetric problem was to actually prove that the circle is the only plane area having this property. A nice proof based on Fourier expansion techniques can be found in [4], pp. 181-187, and a standard proof is offered in [2], pp. 104 ff.

From an elementary viewpoint there is something counterintuitive in the geometrical presentation of this inequality because, on a first and crude approximation, the layman could make a very naïve remark: Draw an area with unit diameter, and then a circle sharing a diameter with it. It "seems obvious" that the whole area is contained in the circle, so the inequality would be an immediate one (see figure 3, left).

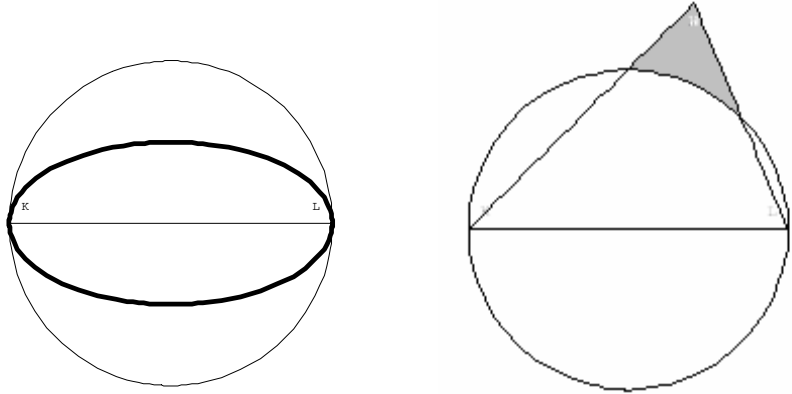

**Figure 3**: The naïve idea and the counterexample

Of course this is erroneous, as is readily shown by considering an isosceles triangle *KLM* with legs *KL* and *KM* longer than the basis *LM*. This is a convex area whose diameter equals one of the legs, say *KL*. Now a circle with *KL* as a diameter is drawn: There is some portion of triangle area in the neighbourhood of corner *M* that lies outside the circle (see figure 3, right)

**A natural question and its answer**

After the above observations, a rather natural question immediately arises: *Is there any upper bound to the fraction of area -of a given convex plane area- lying outside a circle which shares with it a diameter*?

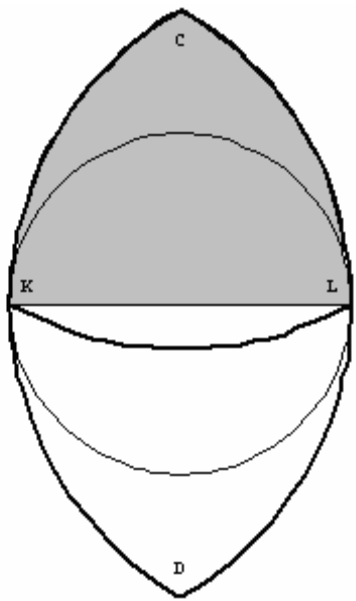

**Figure 4**: Illustrating the procedure.

The classical references [2] and [3] do not mention this topic, and a rather thorough Internet search did not provide direct results, so an attempt to fill this little gap will now be made. In what follows the diameter will be $d = 1$, therefore the area is bounded by $\frac{1}{4}\pi$. To start, consider a unit diameter circle (see figure 4 to follow the discussion), and let *K* and *L* be the endpoints of some diameter thereof.

We shall build a convex area sharing the diameter *KL* with the given circle, and let *X* be any point in the plane. It is evident that both conditions

$$dist(X, K) \leq 1, \; dist(X, L) \leq 1$$

must be satisfied in order that *X* be either an interior or a boundary point of the sought area. Therefore the area is enclosed in the figure defined by two intersecting circle arcs centred at *K* and *L* and having unit radius. They define two points *C* and *D* outside the circle and a curvilinear figure *KDLC*. Of course $dist(C, D) > 1$, so *KDLC* cannot be a candidate for solving our problem. Therefore we restrict our attention to vertex *C* (of course, its symmetrical point *D* could be employed as well) and observe that no point in the boundary of the sought area can be more than 1 apart from it. Although it seems rather natural to draw the unit radius circle arc centred at *C* and joining *K* and *L* to obtain a curvilinear triangle *KLC* as a more appropriate candidate, it is clear that the diameter *KL* does a better job, and we claim that the *mixed* triangle $KLC_{mix}$ is the solution to our problem. Incidentally, the curvilinear triangle *KLC* is called "the Reuleaux triangle" (see Appendix).

A measure of how much area lies outside the circle can now be defined: Simply, it is the ratio between the area outside the circle and the total area just constructed.

$$\mu = \frac{\text{area outside the circle}}{\text{total area}} \leq 1.$$

It is an easy task to compute the total area of the mixed triangle $KLC_{mix}$:

$$\text{total area}(KLC_{mix}) = \frac{\pi}{3} - \frac{\sqrt{3}}{4} = \frac{8\pi - 6\sqrt{3}}{24}$$

and the fraction outside the circle is:

$$\text{exterior area}(KLC_{mix}) = \text{total area}(KLC_{mix}) - \frac{\pi}{8} = \frac{5\pi - 6\sqrt{3}}{24}$$

Therefore, the following value is obtained:

$$\mu = \frac{5\pi - 6\sqrt{3}}{8\pi - 6\sqrt{3}} \cong 0.36$$

*i.e.* the maximum fraction of area lying outside the circle amounts to approximately 36% of the total area.

In order to show optimality of this result, let us consider adding some area to $KLC_{mix}$ by modifying the boundary curves. This cannot be done by changing the curved side *CK* (or *LC*) into another convex curve joining both points for this would imply the existence of some boundary point at a distance from *L* (respectively, *K*) larger than 1, thus contradicting the fact that the sought figure must have unit diameter. Some area can be added below the diameter *KL* preserving both convexity and unit diameter, but in this case the denominator in the definition of the measure would increase, thus reducing the value of $\mu$.

The construction also shows that convexity is a necessary condition for the bound to be a valid one. It is enough to observe (see figure 5) that the area outside the circle is a non-convex figure sharing the unit diameter *KL* with the circle, but 100% of it lies outside the circle, (and indeed is less than $\frac{1}{4}\pi$).

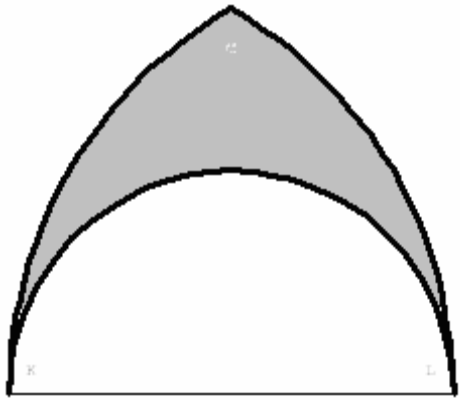

**Figure 5**: Convexity is a necessary condition.

**Appendix on the Reuleaux triangle**

The Reuleaux triangle is a rather familiar curvilinear triangle (see figure A1) obtained by drawing three circle arcs centred at the vertices of an equilateral triangle with a radius equal to the side, it is indeed convex, and its area is the minimum of all possible figures of constant width having the same diameter $d$, a result known as the Blaschke-Lebesgue Theorem. These figures share the common length $\pi \times d$ (a Theorem by Barbier [3]), so the circle and the Reuleaux triangle are extremal curves -in the sense of enclosed area- with this property (see again [3]). Reuleaux triangles have been employed for decorative purposes (see figure A2) and in technological applications, such as the Wankel rotary engines. See also Chapter 7 in [2].

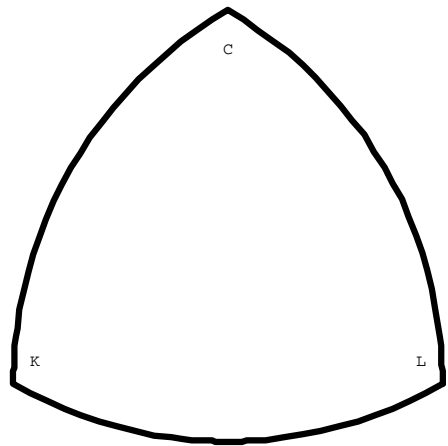

**Figure A1**: The Reuleaux triangle

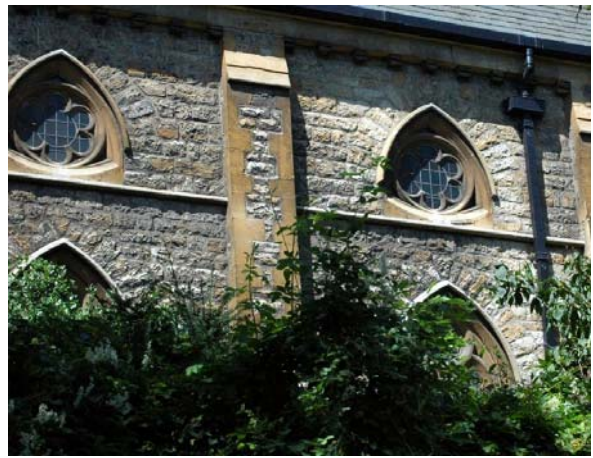

**Figure A2**: Reuleaux triangles in Camden Town, London (photograph by the author).